\documentclass[12pt, reqno]{amsart}
\usepackage{amsmath, amsthm, amscd, amsfonts, amssymb, graphicx,mathrsfs, color}
\usepackage[bookmarksnumbered, colorlinks=true, linkcolor=blue,urlcolor=blue, plainpages]{hyperref}
\input xy \xyoption{all}
\usepackage[all]{xy}
\usepackage{layout}
\usepackage{tikz}
\usepackage{bm}
\newtheorem{theorem}{Theorem}[section]
\newtheorem{lemma}[theorem]{Lemma}
\newtheorem{proposition}[theorem]{Proposition}
\newtheorem{corollary}[theorem]{Corollary}
\theoremstyle{definition}
\newtheorem{definition}[theorem]{Definition}
\newtheorem{example}[theorem]{Example}

\theoremstyle{remark}

\numberwithin{equation}{section}

\begin{document}
\setcounter{page}{1}

\title[polynomial extensions  of $\mathfrak{cP}$-Baer rings]{\textbf{polynomial extensions  of $\mathfrak{cP}$-Baer rings}}

\author[ N. Aramideh  and A. Moussavi  ]{  N. Aramideh  and A. Moussavi$^{*}$  }
\address{Department of Pure Mathematics, Faculty of Mathematical Sciences\\
Tarbiat Modares University, P.O.Box:14115-134, Tehran, Iran}
\email{\textcolor[rgb]{0.00,0.00,0.84}{n.aramideh@modares.ac.ir\,\
aramidehnasibeh@gmail.com}}
\email{\textcolor[rgb]{0.00,0.00,0.84}{ moussavi.a@modares.ac.ir\,\
moussavi.a@gmail.com}}

\thanks{*Corrosponding authour: moussavi.a@modares.ac.ir and moussavi.a@gmail.com} \subjclass{16S34, 16S35, 16S36}
\keywords{right $\mathfrak{cP}$-Baer ring, right p.q.-Baer ring, monoid ring, skew polynomial ring, skew inverse Laurent series ring.}

%_____________________________________________________________________
%_________________________Please type the abstract here_____________________________________________________
\begin{abstract}
 Birkenmeier and Heider, in \cite{Heider}, say that a ring $R$ is right $\mathfrak{cP}$-Baer if the right annihilator of a cyclic projective right $R$-module in $R$ is generated by an idempotent. These rings are a generalization of the right p.q.-Baer and abelian rings. In general, a formal power series ring over one indeterminate, wherein its base ring is right p.q.-Baer, is not necessarily right p.q.-Baer. However, according to \cite{Heider}, if the base ring is right $\mathfrak{cP}$-Baer then the power series ring over one indeterminate is right $\mathfrak{cP}$-Baer. Following \cite{Heider}, we investigate the transfer of the $\mathfrak{cP}$-Baer property between a ring $R$ and various extensions (including skew polynomials, skew Laurent polynomials, skew power series, skew inverse Laurent series, and monoid rings). We also answer a question posed by Birkenmeier and Heider \cite{Heider} and provide examples to illustrate the results. 
\end{abstract}

%%% ----------------------------------------------------------------------
\maketitle
%\tableofcontents
%%% ----------------------------------------------------------------------
%\baselineskip=24pt
\section{Introduction}
Throughout, all rings are associative with an identity. Let $R$ be a ring, $\alpha$ be an  endomorphism and $\delta$ an
$\alpha$-derivation of $R$, that is, $\delta$ be an additive map such that $\delta (ab)=\delta (a)b+\alpha(a)\delta (b)$, for all $a,b\in R$. Recall from \cite{Kaplansky} that  $R$ is  a Baer ring if the right annihilator of every nonempty subset of $R$ is generated by an idempotent. A ring $R$ is called right Rickart (also called right p.p.) if the right annihilator of every  singleton subset of $R$ is generated by an idempotent \cite{Maeda}. Also, from \cite{Clark} and \cite{birkenkimpark}, $R$ is a quasi-Baer  ring if the right annihilator of any ideal is generated by an idempotent as a right ideal. These rings have a rich structure theory because of the abundance of idempotents. Moreover, these rings naturally appear not only in ring theory (e.g., the right Noetherian right semihereditary rings are Baer rings) but also in other areas of Mathematics such as the operators theory (e.g., the von Neumann algebras are Baer rings, and the local multiplier algebras are quasi-Baer rings), see \cite{birkenmeier2000quasi} and \cite{pollingher}. \par
Every right annihilator of a principal ideal is generated by an idempotent in the right p.q.-Baer rings. One of the most important principal right ideals in a ring is one generated by an idempotent. Note that every cyclic projective module is isomorphic to such a right ideal. In \cite{Heider}, Birkenmeier and  Heider define $R$ to be a right $\mathfrak{cP}$-Baer ring if the right annihilator of every cyclic projective $R$-module is generated by an idempotent. They also provide some basic examples and results. For instance, abelian rings (i.e., rings in which every idempotent is central), semisimple, local, biregular, prime, and right p.q.-Baer rings are right $\mathfrak{cP}$-Baer.
In \cite{Heider}, the authors give a large class of rings in their examples which are right $\mathfrak{cP}$-Baer, not left $\mathfrak{cP}$-Baer, and neither right nor left p.q.-Baer. They also provide a large class of rings which are left $\mathfrak{cP}$-Baer, not right $\mathfrak{cP}$-Baer, and neither right nor left p.q.-Baer, and examples of a large class of rings which are right $\mathfrak{cP}$-Baer but not right p.q.-Baer.\par
By \cite{Armendariz}, if $R$ is a Baer and reduced ring, then the polynomial ring $R[x]$ is also Baer. From \cite{birkenkp} and \cite{Zhou}, one can see that the quasi-Baer condition transfers between polynomial rings, formal power series rings, and the base ring with no extra conditions; and the p.q.-Baer condition transfers between polynomial rings and the base ring with no extra conditions. However, the p.q.-Baer condition does not transfer from the base ring to its formal power series extension unless there is an extra condition on the right semicentral idempotents of the base ring \cite{Cheng}. This suggests the question: ``How much" of the right p.q.-Baer condition transfers from the base ring to the formal power series ring? see \cite{Nasr} and\cite{nasr-ore}.\\

Generally, a formal power series ring with one indeterminate whose base ring is right p.q.-Baer is not necessarily right p.q.-Baer. In \cite{Heider}, the authors studied the polynomial and formal power series rings  with respect to the right $\mathfrak{cP}$-Baer condition. Theorem 4.6 in \cite{Heider} states that if the base ring is right $\mathfrak{cP}$-Baer, then the polynomial ring and the formal power series ring over one indeterminate are right $\mathfrak{cP}$-Baer. They raised an open question: Let $R$ be a ring and $X$ a nonempty set of not necessarily commuting indeterminates. If $R$ is
right $\mathfrak{cP}$-Baer, when is $ R[X] $ or $ R[[X]]$ right $\mathfrak{cP}$-Baer?\par
In Subsection \ref{s21}, we  show that for any  $ \alpha $-compatible ring $R$ with an endomorphism $ \alpha $, $R$ is right $\mathfrak{cP}$-Baer if and only if the skew polynomial ring $ R[x;\alpha] $ is right $\mathfrak{cP}$-Baer if and only if
the skew power series ring $ R[[x;\alpha]] $ is right $\mathfrak{cP}$-Baer.\par
In Subsection \ref{s23}, we show that a ring $ R $ is  right $\mathfrak{cP}$-Baer  if and only if the polynomial ring $ R[X] $ is a right $\mathfrak{cP}$-Baer ring if and only if the power series ring $ R[[X]]$  is a right $\mathfrak{cP}$-Baer ring when $ X $ is a nonempty set of not necessarily commuting indeterminates. This result provides the answer to the question posed by Birkenmeier and  Heider in \cite{Heider}.\par
In Subsection \ref{s25}, we investigate the behavior of right $\mathfrak{cP}$-Baer rings with respect to monoid rings.\par
In Subsection \ref{s22}, we provide several examples of skew matrix rings which are $ \alpha $-compatible reversible and right $\mathfrak{cP}$-Baer. We then  show that  an $ \alpha$-compatible  semicommutative ring $R$  is right $\mathfrak{cP}$-Baer  if and only if the skew Laurent polynomial ring $R[x, x^{-1};\alpha] $ is  right $\mathfrak{cP}$-Baer if and only if  Laurent polynomial ring $R[x, x^{-1}] $ is  right $\mathfrak{cP}$-Baer.\par
In Subsection \ref{s3}, we show that for any $ (\alpha,\delta) $-compatible ring  $R$ with an automorphism $ \alpha $ and an  $ \alpha $-derivation $ \delta $, $R$ is right $\mathfrak{cP}$-Baer if and only if the skew inverse power series ring $ R[[x^{-1}; \alpha, \delta]] $ is right $\mathfrak{cP}$-Baer.\par
Examples to demonstrate the theory are provided.
\section{$\mathfrak{cP}$-Baer polynomial  rings }\label{s2}
According to  \cite{Heider}, the class of $\mathfrak{cP}$-Baer rings is defined in terms of principal ideals generated by idempotents. Observe that if $P$ is a cyclic projective right $R$-module, then $P \cong eR$, where $e = e^{2} \in  R$. Hence, $r_{R}(P) = r_{R}(eR) = r_{R}(ReR) = r_{R}(tr(P))$, where $r_{R}(-)$ and $tr(-)$ denote the right annihilator and the trace of $(-)$ in $R$, respectively. A ring $R$ is a right $\mathfrak{cP}$-Baer ring if for every cyclic projective right $ R $-module $P$, $r_{R}(P) = cR$, for some $c = c^{2} \in R$. Therefore, a ring $R$ is right $\mathfrak{cP}$-Baer if and only if for $e = e^{2} \in  R$, $r_{R}(ReR) = r_{R}(eR) = cR$ where  $c = c^{2} \in  R$ if and only if for each $e = e^{2} \in R$ there exists $f \in S_{r}(R)$ such that $r_{R}(eR) = r_{R}(f)$. In addition, we use $I(R)$ and $\mathcal{P}(R)$ to denote the set of idempotents and the prime radical of the ring $R$, respectively.
\subsection{$\mathfrak{cP}$-Baer skew power series  rings }\label{s21}
We denote $ R[x;\alpha] $ the skew polynomial ring whose elements are the polynomials over $ R $, the addition is defined as usual, and the multiplication given by $ xa = \alpha(a)x $, for any $ a\in R $. We also denote to  $R[[x;\alpha]]$, the skew power series ring, whose elements are the series  $ \sum^{\infty}_{i=0} a_{i}x^{i} $, where $ a_{i}\in R $, addition is defined as usual and the multiplication subject to the relation $ xa = \alpha (a)x $ for any $ a\in R $.
\begin{proposition}\cite[Proposition $1.4$]{Heider}\label{6}
	The following conditions are equivalent:
	\begin{enumerate}
		\item R is right $\mathfrak{cP}$-Baer.
		\item If $ e = e^{2}\in R $, then $ r_{R}(ReR) = r_{R}(eR) = cR $ where $ c = c^{2}\in R. $
		\item For each $ e = e^{2}\in R $ there exists $ f = f^{2}\in S_{r}(R) $ such that $ r_{R}(eR) = r_{R}(f). $
		\item $ r_{R}( \sum^{m}_{i=0}e_{i}R) = fR  $, for some $ f = f^{2}\in R $, where $ e_{i} = e^{2}_{i}\in R $ for $ 1\leq i\leq m $.
	\end{enumerate}
\end{proposition}
We begin with a lemma about the behavior of the coefficients of an idempotent skew polynomial and an idempotent skew power series. This lemma and several others are crucial in proving that the right $\mathfrak{cP}$-Baer condition transfers between the base ring and its skew polynomial and skew power series extensions. The following lemma is an extension of \cite[Lemma $4.1$]{Heider}.
\begin{lemma}\label{2}
	Let R be a ring with an endomorphism $ \alpha $ such that $ \alpha (e) = e $ for any $ e = e^{2}\in R $. We have the following statements.
	\begin{enumerate}
		\item If $ e(x) = e_{0} + e_{1}x + \cdots +e_{m}x^{m} $ is an idempotent in $ R[x;\alpha], $ then $ e_{i} \in Re_{0}R $ for any  $ i\geq 0 $.
		\item If $ e(x) = \sum^{\infty}_{i=0} e_{i}x^{i} $ is an idempotent in $ R[[x;\alpha]], $ then $ e_{i} \in Re_{0}R $ for any  $ i\geq 0. $
	\end{enumerate}
\end{lemma}
\begin{proof}
	The proof of (ii) is similar to (i), and we proceed by proving (i). The hypothesis implies that $e_{i} =  \sum^{i}_{k=0} e_{k}\alpha^{k}(e_{i-k}) $ for each $ i\geq 0 $. Hence, $ e_{0} $ is an idempotent of $R$, and so $ \alpha (e_{0}) = e_{0} $. By the induction on $ i $, we show that $ e_{i} \in Re_{0}R $ for any  $ i\geq 0 $. Clearly, $ e_{0}\in Re_{0}R $. As $ e_{1} = e_{0}e_{1} + e_{1}e_{0} $, we have $ e_{1}\in Re_{0}R $. Now we assume that the result holds for all $ i < m $. Then:
	\begin{equation*}
		e_{m} = e_{0}e_{m} + e_{1}\alpha(e_{m-1}) + \cdots + e_{m-1}\alpha^{m-1}(e_{1}) + e_{m}e_{0}.
	\end{equation*}
	Since $ e_{m}e_{0}=e_{m}e_{0}e_{0}\in Re_{0}R $, $ e_{0}e_{m}=e_{0}e_{0}e_{m}\in Re_{0}R  $ and $ e_{i} \in Re_{0}R $ for all  $ 0\leq i\leq m-1 $ so we have, $ e_{m}\in Re_{0}R $. This completes the proof of the lemma.
\end{proof}

An idempotent $ e\in R $ is called left (resp. right) semicentral if $ re = ere$ (resp. $er = ere$), for all $ r\in R $. The set of all left (right)  semicentral idempotents of $R$ are denoted by $ S_{l}(R)$ (resp. $S_{r}(R)$). Define $B(R)= S_{l}(R) \cap S_{r}(R)$ (the set of all central idempotents) and if $R$ is a semiprime ring then $ S_{l}(R) = B(R) = S_{r}(R) $.

\begin{proposition}\label{440}
	Let R be a ring with an endomorphism $ \alpha $ such that $ \alpha (c) = c $ for any left semicentral idempotent $ c \in R $. If $ e(x)\in S_{l}(R[x;\alpha]), $ where $ e(x) = e_{0} + e_{1}x + \cdots +e_{m}x^{m}, $ then:
	\begin{enumerate}
		\item $ e_{0}\in S_{l}(R) $,
		\item $ e_{0}e_{i} = e_{i} $ for all $ 0\leq i\leq m $,
		\item $ e_{i}e_{0} = 0 $ for all $ 1\leq i\leq m $,
		\item $ e(x)R[x;\alpha] = e_{0}R[x;\alpha] $.
	\end{enumerate}
\end{proposition}
\begin{proof}
	If $ m =0 $ we are done. Therefore, assume that $ m\geq 1 $. Since $ e(x)\in S_{l}(R[x;\alpha]) $, $ re(x) = e(x) r e(x) $ for any $ r\in R $. Then we have $ \sum^{m}_{i=0}re_{i}x^{i} =  \sum^{2m}_{k=0}( \sum_{i+j=k}e_{i}\alpha^{i}(re_{j}))x^{k} $. If $ k = 0 $, then $ re_{0} = e_{0}re_{0} $. Thus $ e_{0}\in S_{l}(R) $ and so part (i) is satisfied. Now, we will prove (ii) and (iii) by the induction on $ i $. We know that $ re_{1} = e_{0}re_{1} + e_{1}\alpha (r)e_{0} $. Multiplying on the right by $ e_{0} $, we get $ re_{1}e_{0} = e_{0}re_{1}e_{0} + e_{1}\alpha (r)e_{0} $. As $ re_{1}e_{0} = e_{0}re_{1}e_{0} $,  $ e_{1}\alpha (r)e_{0} = 0 $. Hence, $ re_{1} = e_{0}re_{1} $. Set $ r = 1 $ then $ e_{0}e_{1} = e_{1} $ and $ e_{1}e_{0} = 0 $. Now assume that $ l $ is a positive integer such that $ 1< l\leq m $, and $ e_{i}e_{0} = 0 $ and $ e_{0}e_{i} = e_{i} $ for all $ 1\leq i< l $. Multiplying the equation $ re_{l} = e_{l}\alpha^{l}(r)e_{0} + e_{l-1}\alpha^{l-1}(re_{1})+ \cdots + e_{0}re_{l} $ from the right side by $ e_{0} $, we get $ re_{l}e_{0} = e_{l}\alpha^{l}(r)e_{0} + e_{l-1}\alpha^{l-1}(re_{1})e_{0}+ \cdots + e_{0}re_{l}e_{0} $. Since $ re_{l}e_{0} = e_{0}re_{l}e_{0} $,  $ e_{l}\alpha^{l}(r)e_{0} = 0 $. Thus, $ re_{l} = e_{l-1}\alpha^{l-1}(re_{1})+ \cdots + e_{0}re_{l} $. Multiplying from the left side by $ e_{0} $, we obtain $ e_{0}re_{l} = e_{0}e_{l-1}\alpha^{l-1}(re_{1})+ \cdots + e_{0}re_{l} $. Then $ e_{0}e_{l-1}\alpha^{l-1}(re_{1})+ \cdots + e_{0}e_{1}\alpha(re_{l-1}) = 0 $ and so $ re_{l} = e_{0}re_{l} $. Taking $ r = 1 $ yields $ e_{0}e_{l} = e_{l} $ and $ e_{l}e_{0} = 0 $. On the other hand, by  (ii) and (iii), $ e_{0}e(x) = e(x) $ and $ e(x)e_{0} = e_{0} $. Therefore, $ e(x) \in e_{0}R[x;\alpha] $ and $ e_{0}\in e(x)R[x;\alpha] $, so $ e(x)R[x;\alpha] = e_{0}R[x;\alpha] $ and we are done.
\end{proof}
For the case of skew power series rings, an argument similar to the above yields the following.
\begin{proposition}\label{4}
	Let R be a ring with an endomorphism $ \alpha $ such that $ \alpha (c) = c $ for any left semicentral idempotent $ c \in R $. If $ e(x)\in S_{l}(R[[x;\alpha]]), $ where $ e(x) = \sum^{\infty}_{i=0} e_{i}x^{i}, $ then:
	\begin{enumerate}
		\item $ e_{0}\in S_{l}(R) $,
		\item $ e_{0}e_{i} = e_{i} $ for all $ i\geq 0 $,
		\item $ e_{i}e_{0} = 0 $ for all $ i\geq 1 $,
		\item $ e(x)R[[x;\alpha]] = e_{0}R[[x;\alpha]] $.
	\end{enumerate}
\end{proposition}

\begin{lemma}\label{7}
	Let R be a right $\mathfrak{cP}$-Baer ring with a monomorphism $ \alpha $ such that $ \alpha (c) = c $ for any left semicentral idempotent $ c\in R $. If $ eR\alpha^{k}(a) = 0 $ for some integer $ k\geq 0 $, where $ e = e^{2}\in R $ and $ a\in R $, then $ eR\alpha^{m}(a) = 0 $ for any integer $ m\geq 0 $.
\end{lemma}
\begin{proof}
	Let $ eR\alpha^{k}(a) = 0 $ for some integer $ k\geq 0 $, where $ e = e^{2}\in R $ and $ a\in R $. By the assumption, there exists $ c = c^{2}\in R $ such that $ r_{R}(eR) = cR $. So $ \alpha^{k}(a)\in cR $. Note that since $ r_{R}(eR) $ is an ideal in $ R $, $ c\in S_{l}(R) $. Hence, $ \alpha^{k}(a) = c\alpha^{k}(a) $ and $ a = ca $. Therefore, $ eR\alpha^{m}(a) = eR\alpha^{m}(ca) = eRc\alpha^{m}(a) = 0, $ for any integer $ m\geq 0 $.
\end{proof}

\begin{lemma}\label{8}
	Let R be a right $\mathfrak{cP}$-Baer ring with a monomorphism $ \alpha $ such that $ \alpha (c) = c $ for any left semicentral idempotent $ c\in R $. For any idempotent $ e(x) = \sum^{m}_{i=0}e_{i}x^{i}\in  R[x;\alpha]$ and any $ p(x) = \sum^{n}_{j=0}a_{j}x^{j}\in  R[x;\alpha] $, if $ e(x)R[x; \alpha]p(x) = 0 $ then $ e_{0}R\alpha^{l}(a_{j}) = 0 $ for any integer $ l\geq 0 $ and $ 0\leq j\leq n $.
\end{lemma}
\begin{proof}
	Suppose that $ e(x)R[x;\alpha]p(x) = 0 $, where $e(x) = \sum^{m}_{i=0}e_{i}x^{i}\in  R[x;\alpha]$ is an idempotent and $ p(x) = \sum^{n}_{j=0}a_{j}x^{j}\in  R[x;\alpha] $. So, $ e(x)rp(x) = 0 $ for all $ r\in R $. Hence, $ 0 = \sum^{m+n}_{k=0}( \sum_{i+j=k}e_{i}\alpha^{i}(ra_{j}))x^{k} = c_{0} + \cdots + c_{m+n}x^{m+n} $. First, we will prove that $ e_{0}Ra_{j} = 0 $ for each $ 0\leq j\leq n $. We utilize the induction on $ j $. If $ j = 0 $, then $ 0 = c_{0} = e_{0}ra_{0} $. Hence, $ e_{0}Ra_{0} = 0 $. Now assume the result is true for all $1\leq t < j $. So $ e_{0}Ra_{t} = 0 $ for each $1\leq t < j $. Hence, by Lemma \ref{7}, $ e_{0}R\alpha^{l}(a_{t}) = 0 $ for any integer $ l\geq 0 $ and each $1\leq t < j $. We have
	\begin{equation}\label{e1}
		0= c_{j} = e_{0}ra_{j} + e_{1}\alpha (ra_{j - 1})+ \cdots +e_{j}\alpha^{j}(ra_{0}).
	\end{equation}
	Multiplying $ e_{0} $ from the left side of (\ref{e1}), we obtain
	\begin{equation*}
		0 = e_{0}ra_{j} + e_{0}e_{1}\alpha (r)\alpha(a_{j - 1})+ \cdots +e_{0}e_{j}\alpha^{j}(r)\alpha^{j}(a_{0}).
	\end{equation*}
	Now, the induction hypothesis gives $ e_{0}Ra_{j} = 0$. Therefore, by Lemma \ref{7}, $ e_{0}R\alpha^{l}(a_{j}) = 0 $ for any integer $ l\geq 0 $ and $ 0\leq j\leq n $.
\end{proof}
\begin{lemma}\label{88}
	Let R be a right $\mathfrak{cP}$-Baer ring with a monomorphism $ \alpha $ such that $ \alpha (c) = c $ for any left semicentral idempotent $ c\in R $. For any idempotent $ e(x) = \sum^{\infty}_{i=0}e_{i}x^{i}\in  R[[x;\alpha]]$ and any $ p(x) = \sum^{\infty}_{j=0}a_{j}x^{j}\in  R[[x;\alpha]] $, if $ e(x)R[[x;\alpha]]p(x) = 0 $ then $ e_{0}R\alpha^{l}(a_{j}) = 0 $ for any integer $ l\geq 0 $ and $ j\geq 0 $.
\end{lemma} 
\begin{proof}
	The proof is similar to the proof of Lemma \ref{8}.
\end{proof}
According to Krempa \cite{Krempa}, an endomorphism $\alpha$ of a ring $R$ is said to be \textit{rigid} if $a\alpha (a) = 0$ implies $a = 0$, for $a \in R$. A ring $R$ is said to be $\alpha$-\textit{rigid} if there exists a rigid endomorphism $\alpha$ of $R$. Clearly, every domain $D$ with a monomorphism $\alpha$ is rigid. In \cite{Hashemi}, the authors introduced $(\alpha,\delta)$-compatible rings and studied their properties. A ring $R$ is $\alpha$-\textit{compatible} if for each $a,b \in R$, $ab = 0$ if and only if $a \alpha (b) = 0$. This may only happen when the endomorphism $\alpha$ is injective. Moreover, $R$ is said to be $\delta$-\textit{compatible} if for each $a,b \in R$, $ab = 0 \Rightarrow a\delta(b) = 0$. A ring $R$ is $(\alpha,\delta)$-\textit{compatible} if it is both $\alpha$-compatible and $\delta$-compatible. Also, by \cite[Lemma $ 2.2 $]{Hashemi}, a ring $R$ is $\alpha$-rigid if and only if $R$ is $(\alpha,\delta)$-compatible and reduced.
By \cite{Hashemi}, whenever $R$ is an $ (\alpha,\delta) $-compatible ring, then $ \alpha (e) = e$ and $ \delta (e) = 0$, for any  idempotent $ e\in R $. By \cite{Hashemi}, when $ \alpha $ is a compatible endomorphism of a ring $R$, then $ \alpha $ is injective.
Also, $ \alpha $ is compatible if and only if for $ a,b\in R $, $ \alpha(a)b = 0 \Leftrightarrow ab = 0 $.\par

By \cite{Hashemi}\label{11}, when  $R$ is a ring with an endomorphism $ \alpha $ such that for any $ a,b\in R $, $ ab = 0 $ implies $ a\alpha (b) = 0,$ then $ \alpha (e) = e $ for every idempotent $ e\in R $, see also \cite{nasr-weakly rigid}.

\begin{theorem}\label{12}
	Let R be a ring with a compatible  endomorphism $ \alpha $. Then the following statements are equivalent: 
	\begin{enumerate}
		\item R is right $\mathfrak{cP}$-Baer.
		\item $ R[x;\alpha] $ is right $\mathfrak{cP}$-Baer.
		\item $ R[[x;\alpha]] $ is right $\mathfrak{cP}$-Baer.
	\end{enumerate}
\end{theorem}
\begin{proof}
	(i) $ \Rightarrow $ (ii)  Let $R$ be a right $\mathfrak{cP}$-Baer ring and $ e(x) = \sum^{m}_{i=0}e_{i}x^{i}\in  R[x;\alpha]$ be an idempotent. Then $ e_{0} = e^{2}_{0}\in R $ and so, there exists $ c = c^{2}\in R $ such that $ r_{R}(e_{0}R) = cR $. Since $ r_{R}(e_{0}R) $ is an ideal in $ R $, $  c\in S_{l}(R) $. We claim that $ r_{R[x;\alpha]}(e(x)R[x;\alpha]) = cR[x;\alpha] $. Assume $ p(x) = \sum^{n}_{j=0}a_{j}x^{j}\in r_{R[x;\alpha]}(e(x)R[x;\alpha]) $. Then $ e(x)R[x;\alpha]p(x) = 0  $. So, by Lemma \ref{8}, $ e_{0}Ra_{j} = 0 $ for all $ 0\leq j\leq n $.  Hence, $ a_{j}\in cR $ for all $ 0\leq j\leq n $ and so $ a_{j} = ca_{j} $. Therefore, $ cp(x) = \sum^{n}_{j=0}ca_{j}x^{j} = p(x) $. Then $ p(x)\in cR[x;\alpha] $ and so $ r_{R[x;\alpha]}(e(x)R[x;\alpha])\subseteq  cR[x;\alpha]$. Conversely, it is enough to show that $ c\in r_{R[x;\alpha]}(e(x)R[x;\alpha]) $. Let $ p(x) = \sum^{n}_{j=0}a_{j}x^{j}\in  R[x;\alpha] $ be an arbitrary element. So,
	\begin{align*}
		e(x)p(x)c &= (\sum^{m}_{i=0}e_{i}x^{i})(\sum^{n}_{j=0}a_{j}x^{j})c\\
		&= \sum^{m+n}_{k=0}( \sum_{i+j=k}e_{i}\alpha^{i}(a_{j})c)x^{k}.
	\end{align*}
	By Lemma \ref{2}, $ e_{i}\in Re_{0}R $ for all $ 0\leq i\leq m $. Therefore, $ e_{i}a_{j}c = 0 $ for all $ 0\leq i\leq m $ and $ 0\leq j\leq n $. Thus, $ c\in r_{R[x;\alpha]}(e(x)R[x;\alpha]) $, which shows that $ cR[x;\alpha]\subseteq r_{R[x;\alpha]}(e(x)R[x;\alpha]) $.\\
	(ii) $ \Rightarrow $ (i)  Assume that $ R[x;\alpha] $ is a right $\mathfrak{cP}$-Baer ring and $ c = c^{2}\in R $. Then $r_{R[x;\alpha]}(cR[x;\alpha]) = e(x) R[x;\alpha]$, for an idempotent $ e(x) = \sum^{m}_{i=0}e_{i}x^{i}\in R[x;\alpha]. $ By Proposition \ref{440}, we have $ e(x) R[x;\alpha] = e_{0}R[x;\alpha] $ where $ e_{0}\in S_{l}(R) $. Thus, $ r_{R[x;\alpha]}(cR[x;\alpha]) = e_{0}R[x;\alpha] $. Therefore, $ r_{R}(cR) = e_{0}R $ and so $R$ is a right $\mathfrak{cP}$-Baer ring.\\
	(i) $ \Rightarrow $ (iii) Let $ e(x) = \sum^{\infty}_{i=0}e_{i}x^{i}\in  R[[x;\alpha]] $ be an idempotent. Then $ e_{0} = e^{2}_{0}\in R $ and so, there exists an idempotent $ c\in R $ such that $ r_{R}(e_{0}R) = cR $. We will show that $ r_{R[[x;\alpha]]}(e(x)R[[x;\alpha]]) = cR[[x;\alpha]] $. Using part (ii) of Lemma \ref{2}, $ cR[[x;\alpha]]\subseteq r_{R[[x;\alpha]]}(e(x)R[[x;\alpha]]) $.  For the converse,  we know that
	 $ r_{R[[x;\alpha]]}(e(x)R[[x;\alpha]])\subseteq r_{R[[x;\alpha]]}(e(x)R), $ so it is sufficient to show that $ r_{R[[x;\alpha]]}(e(x)R)\subseteq  cR[[x;\alpha]]$. Let $ p(x) = \sum^{\infty}_{j=0}a_{j}x^{j}\in r_{R[[x;\alpha]]}(e(x)R) $. Then we have $ e(x)Rp(x) = 0 $, and so $ e_{0}Ra_{j} = 0 $ for any $ j\geq 0 $, by Lemma \ref{88}. Hence, $ a_{j}\in cR $ for each $ j\geq 0 $. Thus, $ p(x)\in cR[[x;\alpha]] $. Therefore, $ r_{R[[x;\alpha]]}(e(x)R[[x;\alpha]]) = cR[[x;\alpha]] $ and so  $ R[[x;\alpha]] $ is right $\mathfrak{cP}$-Baer.	\\
 (iii) $ \Rightarrow $ (i) It can be proved by utilizing Proposition \ref{4} and employing a similar procedure to (ii) $ \Rightarrow $ (i).
\end{proof}

\begin{corollary}\label{1o}\cite[Theorem 4.6.]{Heider}
	The following statements are equivalent: 
	\begin{enumerate}
		\item R is a right $\mathfrak{cP}$-Baer ring.
		\item $ R[x] $ is right $\mathfrak{cP}$-Baer.
		\item $ R[[x]] $ is right $\mathfrak{cP}$-Baer.
	\end{enumerate}
\end{corollary}

From \cite[Example 3.6]{birkenkimpark}, there exists a ring $R$ that is p.q.-Baer (and p.p.) while the formal power series ring $ R[[x]] $ is not p.q.-Baer (and not p.p.).

\begin{corollary}\label{1oo}
	Let R be a  right p.q.-Baer ring with a compatible automorphism $ \alpha$. Then $ R[[x;\alpha]] $ is a right $\mathfrak{cP}$-Baer ring.
\end{corollary}
\begin{corollary}\label{1ooo}\cite[Corollary 4.7.]{Heider}
	If R is a  right p.q.-Baer ring, then $ R[[x]] $ is a right $\mathfrak{cP}$-Baer ring.
\end{corollary}

The following example shows that Theorem \ref{12} is not true in general without the compatibility condition.
\begin{example}\cite[Example $ 1 $]{h1}\label{13}
	Let $ F $ be a field and $ R = F1 + \oplus_{i\in \mathbb{Z}}F_{i} $ where $ F_{i} = F $ for all $ i\in \mathbb{Z} $. Thus, $R$ is a commutative regular ring and so $R$ is a $\mathfrak{cP}$-Baer ring. Furthermore, assume that $ \alpha: R\longrightarrow R $ is an automorphism defined by $ \alpha ((a_{i})) = (a_{i+1}) $ for each $i\in \mathbb{Z}$. Let $ a = (a_{i}) $ where $ a_{1} = 1 $ and $ a_{i} = 0 $ for any $ i\neq 1 $ and also $ b = (b_{i}) $ where $ b_{0} = 1 $ and $ b_{i} = 0 $ for any $ i\neq 0 $. Then $ ab=0 $ but $ a\alpha (b) = b\neq 0 $. Thus, R is not $ \alpha $-compatible.
	Now we  show that $ R[x;\alpha] $ is not  $\mathfrak{cP}$-Baer. Assume, to the contrary, that $ R[x;\alpha] $ is a $\mathfrak{cP}$-Baer ring. Let  $ c = (c_{i})\in R $ such that $ c_{0} = 1 $ and $ c_{i} = 0 $ for all $ i\neq 0. $ As $ r_{R[x;\alpha]}(cR[x;\alpha]) $ contains $ \oplus_{i\geq 1}F_{i} ,$  $  r_{R[x;\alpha]}(cR[x;\alpha]) \neq 0 $. Thus, there exists $e^{2}(x) = e(x) = \sum^{m}_{i=0}e_{i}x^{i}\in  R[x;\alpha] $ such that $ r_{R[x;\alpha]}(cR[x;\alpha]) = e(x)R[x;\alpha] $. Hence, $ 0 = crx^{t}e(x) = \sum^{m}_{i=0}cr\alpha^{t}(e_{i})x^{i+ t} $ for any $ r\in R $ and any integer $ t\geq 0 $. Thus, $ cr\alpha^{t}(e_{m}) = 0 $ for any integer $ t\geq 0 $. So $ e_{m} = (e_{m_{j}}) $ where $ e_{m_{j}} = 0 $ for all $ j\leq 0 $ and $ e_{m_{j}} $ is eventually $ 0 $. We note that if $ r = 1 $ then $ c\alpha^{k}(e_{m}) = 0 $ for any integer $ k\geq 0 $. So $ 0 = crx^{t}e(x) = \sum^{m - 1}_{i=0}cr\alpha^{t}(e_{i})x^{i+ t} $. Continuing this process, we have $ e_{i} = (e_{i_{j}})\in R $ for all $ 0\leq i\leq m $, where $ e_{i_{j}} = 0 $ for all $ j\leq 0 $ and there exists $ t_{i} $ such that $ e_{i_{t_{k}}} = 0 $ for all $ t_{i}\leq t_{k} $. Let $ t = \max\{t_{0},\cdots,t_{m}\} $ and set $ h = (h_{l}) $ where $ h_{2t} = 1 $ and $ h_{l} = 0 $ for any $ l\neq 2t $. Then $ h\in r_{R[x;\alpha]}(cR[x;\alpha]) = e(x)R[x;\alpha] $. So there exists $ g(x) = \sum^{n}_{j=0}g_{j}x^{j}\in R[x;\alpha] $ such that $ h = e(x)g(x) $. Therefore, $ h = e_{0}g_{0} $ and so the $ (2t)th $ component of $ h $ is zero, which is a contradiction.
\end{example}
From \cite [Definition $ 2.1 $]{nasr2011skew}\label{130}, a ring $R$ with an endomorphism $ \alpha $ is called skew powerserieswise Armendariz $(or\ SPA\ for\ short)$, if for every skew power series $ f(x) = \sum^{\infty}_{i=0}a_{i}x^{i} $ and $ g(x) = \sum^{\infty}_{j=0}b_{j}x^{j} $ in $ R[[x;\alpha]] $, $ f(x)g(x) = 0 $ if and only if $ a_{i}b_{j} = 0 $ for all $ i,j $.
By \cite [Lemma $ 2.2 $]{nasr2011skew}\label{131}, every $SPA$-ring $R$ is $  \alpha $-compatible for every endomorphism $  \alpha $ of $R$. Using Theorem \ref{12} and Lemma \ref{2}, we get the following.

\begin{corollary}\label{132}
	Let R be an SPA-ring with an endomorphism $ \alpha $.  Then the following statements are equivalent: 
	\begin{enumerate}
		\item R is right $\mathfrak{cP}$-Baer ring.
		\item $ R[x;\alpha] $ is right $\mathfrak{cP}$-Baer ring.
		\item $ R[[x;\alpha]] $ is right $\mathfrak{cP}$-Baer ring.
	\end{enumerate}
\end{corollary}
\subsection{$\mathfrak{cP}$-Baer polynomial rings with several variables}\label{s23}

From \cite{birk paa}, it is shown that the quasi-Baer condition transfers between $R$, $R[X]$ and $R[[X]]$
where $X$ is a nonempty set of not necessarily commuting indeterminates. In \cite{birkenmeier2003triangular}, we see that the p.q.-Baer condition transfers to $R[X]$, but in general, it does not transfer to $R[[X]]$. \par
According to  \cite[Theorem 4.6]{Heider}, if the base ring is right $\mathfrak{cP}$-Baer, then the polynomial ring and the formal power series ring over one indeterminate is right $\mathfrak{cP}$-Baer. In  \cite{Heider},  Birkenmeier and  Heider  raised an open question: when $ R[X] $ and $ R[[X]]$ are right $\mathfrak{cP}$-Baer if the base ring is right $\mathfrak{cP}$-Baer, where $X$ is a nonempty set of not necessarily commuting indeterminates?

\begin{lemma}\label{40}
	Let R be a ring and $ X $ a nonempty set of not necessarily commuting indeterminates.  If $ e $ is an idempotent in $ R[X] $ or $ R[[X]] $ with the constant term $ e_{0}\in R $, then coefficients of $ e $ all belong to $ Re_{0}R $.
\end{lemma}
\begin{proof}
	 Let $ e = e^{2} $, when $ X = \lbrace x\rbrace $. Then  the result follows immediately from Lemma \ref{2}. We will prove the result for $ X = \lbrace x_{1},x_{2}\rbrace $ and $ R[[X]] $. The proof for an arbitrary nonempty set of indeterminates $ X $ and $ R[X] $ is a straightforward generalization. Assume that $ e =  \sum_{i = 0}^{\infty}\sum^{\infty}_{j = 0}e_{(i,j)}x^{i}_{1}x^{j}_{2} $ with $ e_{0} = e_{(0,0)} $ the constant term of $ e $. We claim that $ e_{(i,j)} \in Re_{0}R $ for all $ i,j \geq 0 $, and we prove this by induction. Since $ e^{2} = e $, $ e^{2}_{0} = e_{0} $ and so $ e_{0}\in Re_{0}R $. Now,
	 \begin{equation*}
	  e_{0}e_{(1,0)} + e_{(1,0)}e_{0} = e_{(1,0)}; \quad   e_{0}e_{(0,1)} + e_{(0,1)}e_{0} = e_{(0,1)}.
	 \end{equation*}
 So $ e_{(1,0)}$ and $ e_{(0,1)} $ belong to $ Re_{0}R $. Assume that the result holds for all $ 0 \leq k < i $ and $ 0 \leq l < j $. Note that
 \begin{equation*}
 	 e_{0}e_{(i,j)} + e_{(i,0)}e_{(0,j)} + e_{(i,1)}e_{(0,j-1)} + \cdots + e_{(0,j)}e_{(i,0)}+e_{(i,j)}e_{0} = e_{(i,j)}.
 \end{equation*}
 Since ,by the induction hypothesis, $ e_{(k,l)} \in Re_{0}R $ for all $ 0 \leq k < i $ and $ 0 \leq l < j $ so, $ e_{(i,j)}\in Re_{0}R $. Therefore, the proof is complete.
\end{proof}
 We recall the following result of Birkenmeier et,al., which will be used in our main result in this subsection.
\begin{theorem}\cite[Theorem  2.3]{birkenmeier2003triangular}\label{36}
	Let $ R $ be a ring and $ X $ a nonempty set of not necessarily commuting indeterminates. Let $ \Delta = R[X] $ or $ R[[X]] $. Then $ e \in S_{l}(\Delta) $ with $ e_{0} \in R $ the constant term of $ e $ if and only if $ e_{0} \in S_{l}(R) $, $ e_{0}e = e $, and $ e_{0} = ee_{0 }$.
\end{theorem}

\begin{theorem}\label{38}
	Let $ R $ be a ring and $ X $ a nonempty set of not necessarily commuting indeterminates. Then the following conditions are equivalent:
	\begin{enumerate}
		\item $ R $ is a right $\mathfrak{cP}$-Baer ring.
		\item $ R[X] $ is a right $\mathfrak{cP}$-Baer ring.
		\item $ R[[X]] $ is a right $\mathfrak{cP}$-Baer ring.
	\end{enumerate}
\end{theorem}
\begin{proof}
	When $ X = \lbrace x\rbrace $, the result is true by Corollary \ref{1o}. The results will be demonstrated for $ X = \lbrace x_{1},x_{2}\rbrace $. The proof for an arbitrary nonempty set of indeterminates $ X $ is a straightforward generalization.\\
	(i) $ \Rightarrow $ (ii) 
	Let us suppose that $ e =  \sum^{m}_{i = 0}\sum^{n}_{j = 0}e_{(i,j)}x^{i}_{1}x^{j}_{2}  $ is an idempotent of $ R[X] $ and that $ e_{0} = e_{(0,0)} $ is the constant term of $ e $. Since $ e_{0} $ is an idempotent in $ R $, there exists $ c = c^{2}\in R $ such that $ r_{R}(e_{0}R) = cR $. We will show that $ r_{R[X]}(eR[X]) = cR[X] $. First, we show that $ cR[X] \subseteq r_{R[X]}(eR[X]) $. For this purpose, it is enough to show $ c\in r_{R[X]}(eR[X]) $. Let $ f =  \sum^{p}_{k = 0}\sum^{q}_{l = 0}a_{(k,l)}x^{k}_{1}x^{l}_{2} $ be an arbitrary element of $ R[X] $. We know that 
	$$ efc = (\sum^{m}_{i = 0}\sum^{n}_{j = 0}e_{(i,j)}x^{i}_{1}x^{j}_{2})(\sum^{p}_{k = 0}\sum^{q}_{l = 0}a_{(k,l)}x^{k}_{1}x^{l}_{2})c = \sum^{m}_{i = 0}\sum^{n}_{j = 0}\sum^{p}_{k = 0}\sum^{q}_{l = 0}(e_{(i,j)}a_{(k,l)}c)x^{i}_{1}x^{j }_{2}x^{k}_{1}x^{l}_{2}.$$
	On the other hand, $ e_{(i,j)}Rc = 0 $ for all $ i,j\geq 0 $, by Lemma \ref{40}. Then $ efc = 0 $. Therefore, $ c\in r_{R[X]}(eR[X]) $ and so $ cR[X]\subseteq r_{R[X]}(eR[X]) $.\\
	Note that $ r_{R[X]}(eR[X]) \subseteq r_{R[X]}(eR) $. We claim that $ r_{R[X]}(eR)\subseteq cR[X] $. Let $ f =  \sum^{p}_{k = 0}\sum^{q}_{l = 0}a_{(k,l)}x^{k}_{1}x^{l}_{2} $ be an arbitrary element in $ r_{R[X]}(eR) $. Then $ 0 = erf = \sum^{m}_{i = 0}\sum^{n}_{j = 0}\sum^{p}_{k = 0}\sum^{q}_{l = 0}(e_{(i,j)}ra_{(k,l)})x^{i}_{1}x^{j }_{2}x^{k}_{1}x^{l}_{2} $ for every element $ r $ of $ R $. First we will prove that $ e_{0}ra_{(k,l)} = 0 $ for each $ k,l\geq 0 $. We proceed by induction
	on  $ k, l $. Note  that $ e_{0}ra_{(0,0)} = 0 $.
	Now assume the result is true for all $1\leq s < k $ and $1\leq t < l $. So $ e_{0}ra_{(s,t)} = 0 $ for each $1\leq s < k $ and $1\leq t < l $. We have
	\begin{equation*}\tag{$ \ast $}
		 e_{0}ra_{(k,l)} + e_{(1,0)}ra_{(k - 1,l)}+ \cdots +e_{(k,0)}ra_{(0,l)}+ e_{(1,1)}ra_{(k-1,l-1)} + \cdots e_{(k,l)}ra_{(0,0)} = 0.
	\end{equation*}
    By multiplying $ (\ast) $ on the left side by $ e_{0} $, we get
	\begin{equation*}
	 e_{0}ra_{(k,l)} + e_{0}e_{(1,0)}ra_{(k - 1,l)}+ \cdots +e_{0}e_{(k,0)}ra_{(0,l)}+ e_{0}e_{(1,1)}ra_{(k-1,l-1)} + \cdots e_{0}e_{(k,l)}ra_{(0,0)} = 0.
	\end{equation*}
   But since $ e_{0}ra_{(s,t)} = 0 $ for $1\leq s < k $ and $1\leq t < l $ by the induction hypothesis, it follows that $ e_{0}ra_{(k,l)} = 0$. Therefore, $ e_{0}ra_{(k,l)} = 0 $ for each $ k,l\geq 0 $. Thus, $ a_{(k,l)}\in cR $ and so $ ca_{(k,l)} = a_{(k,l)} $. Hence, $ cf =  \sum^{p}_{k = 0}\sum^{q}_{l = 0}ca_{(k,l)}x^{k}_{1}x^{l}_{2} = \sum^{p}_{k = 0}\sum^{q}_{l = 0}a_{(k,l)}x^{k}_{1}x^{l}_{2} = f $ and so $ f\in cR[X] $. Therefore, $R[X] $ is a right $\mathfrak{cP}$-Baer ring.\\  
   (i) $ \Rightarrow $ (iii) Let $ e $ be an idempotent of $ R[[X]] $ and $ e_{0} $ its constant term. Since $ e_{0}\in R $ is an idempotent, there exists $ c = c^{2}\in R $ such that $ r_{R}(e_{0}R) = cR $. It suffices to show that $ r_{R[[X]]}(eR[[X]]) = cR[[X]] $. First, to see that $ cR[[X]]\subseteq r_{R[[X]]}(eR[[X]]) $, take $ f = \sum^{\infty}_{k = 0}\sum^{\infty}_{l = 0}a_{(k,l)}x^{k}_{1}x^{l}_{2}\in eR[[X]] $. We will show that $ fc =  \sum^{\infty}_{k = 0}\sum^{\infty}_{l = 0}a_{(k,l)}cx^{k}_{1}x^{l}_{2} = 0 $. If $ f = 0 $ then $ fc = 0 $. So we assume that $ f \neq 0 $. If $ a_{(0,0)} \neq 0 $ then $ a_{(0,0)}\in e_{0}R $. Thus, $ a_{(0,0)}c = 0 $. If $ a_{(1,0)}c \neq 0 $ then $ a_{(1,0)}c\in Re_{0}R $, by Lemma \ref{40}. But, $ a_{(1,0)}c = (a_{(1,0)}c)c\in Re_{0}Rc = 0 $, which is a contradiction. Similarly, we can show that $ a_{(i,j)}c = 0 $ for all $ i,j\geq 0 $. Hence, $ fc = 0 $ and so $ c\in r_{R[[X]]}(eR[[X]]) $. Therefore, $ cR[[X]]\subseteq r_{R[[X]]}(eR[[X]]) $.\par
   For the converse, we will prove that $ r_{R[[X]]}(eR[[X]])\subseteq cR[[X]] $. Assume that $ f =  \sum^{\infty}_{k = 0}\sum^{\infty}_{l = 0}a_{(k,l)}x^{k}_{1}x^{l}_{2}\in r_{R[[X]]}(eR[[X]]) $. We shall show that $ f = c f$. Let $ g = \sum^{\infty}_{i = 0}\sum^{\infty}_{j = 0}b_{(i,j)}x^{i}_{1}x^{j}_{2} \in eR[[X]]$. Hence, $ gf = 0 $ and so $ b_{(0,0)}a_{(0,0)} = 0 $. As $ b_{(0,0)}\in e_{0}R $, we have $ a_{(0,0)}\in r_{R}(e_{0}R) = cR\subseteq r_{R[[X]]}(eR[[X]]) $. Thus, $ a_{(0,0)} = ca_{(0,0)} $ and $ b_{(i,j)}a_{(0,0)} = 0 $ for all $ i,j\geq 0 $. On the other hand, $ b_{(1,0)}a_{(0,0)}x_{1} + b_{(0,0)}a_{(1,0)}x_{1} = 0 $. So $ b_{(0,0)}a_{(1,0)} = 0 $. Again since $ a_{(1,0)} \in r_{R}(e_{0}R) = cR\subseteq r_{R[[X]]}(eR[[X]]) $, $ a_{(1,0)} = ca_{(1,0)} $ and $ b_{(i,j)}a_{(1,0)} = 0 $ for all $ i,j\geq 0 $. Continuing in this process, we get $ a_{(k,l)} = ca_{(k,l)} $ for each $ k,l\geq 0 $. Therefore, $ f = c f$, and we are done.\\
	(ii) $ \Rightarrow $ (i) Let $ R[X] $ be a right $\mathfrak{cP}$-Baer ring and $ c = c^{2}\in R $. Thus, $ r_{R[X]}(cR[X]) = eR[X] $ where $ e \in R[X]$ is an idempotent with $ e_{0}\in R $ the constant term of $ e $. Note that since $ r_{R[X]}(cR[X]) $ is an ideal in R, $ e\in S_{l}(R[X]) $. By Theorem \ref{36}, $ e_{0}\in S_{l}(R) $ so it is enough to show that $ r_{R}(cR) = e_{0}R $. Let $ a\in r_{R}(cR)\subseteq r_{R[X]}(cR[X]) $. Thus, $ a\in eR[X] $ and so $ ea = a $. Hence, $ e_{0}a = a $. Therefore, $ a\in e_{0}R $ and so $ r_{R}(cR)\subseteq e_{0}R $. Conversely, we have $ r_{R[X]}(cR[X]) = eR[X] $ so $ e\in r_{R[X]}(cR[X]) $. Thus, $ cRe = 0 $ and so $cRe_{0} = 0 $. Therefore, $ e_{0}R\subseteq r_{R}(cR) $, and that $ R $ is a right $\mathfrak{cP}$-Baer ring.\\
	(iii) $ \Rightarrow $ (i) By Theorem \ref{36} and a similar argument as employed in Theorem \ref{12}(ii) $ \Rightarrow $ (i), the result follows.
\end{proof}
\subsection{$\mathfrak{cP}$-Baer monoid rings}\label{s25}
A monoid $ G $ is said to be ordered if the elements of $ G $ are linearly ordered with respect to the relation $ < $ and that, for all $ g, h, s \in G $, $ g < h $ implies $ sg < sh $ and $ gs < hs $. It is an immediate consequence of this  that if $ g < h $ and $ g^{\prime} < h^{\prime} $, then $ gg^{\prime} < hh^{\prime} $. By using \cite[Lemma 13.1.6 and Corollary 13.2.8]{pass}, it is acknowledged that ordered groups consist of torsion-free nilpotent groups and free groups. Also, we call the ordered monoid $ G $ a positive monoid with an identity $ \mu $ whenever  $ \mu \leq g $ for every $ g\in G $.

\begin{lemma}\label{l1}
	Let $ G $ be a positive monoid with an identity $ \mu $. If $ g,h,k \in G $ such that $ k\neq \mu $ and $ k = gh $, then $ g, h \leq k$.
\end{lemma}
\begin{proof}
	Assume, to the contrary, that it is not. So at least one of $ g $ or $ h $ is greater than $ k $. Consider the following two cases:
	\begin{itemize}
		\item[(1)] If $ k < g , h $, then by hypothesis, $ k = k\mu < kk< gh = k $. Hence, $ k = k^{2} $, and so $ k = \mu $, which is a contradiction.
		\item[(2)] If $ g < k < h $ then, $ gk < gh = k $. Hence, $ k = \mu k < gk < k $, and so $ gk = k $. Therefore, $ g = \mu $, and so $ k = h $, which is a contradiction.
	\end{itemize}
\end{proof}

Furthermore, a monoid $G$ is said to be cancellative if for all $ g,h\in G $, $ sg = sh\ or\ gs = hs $ for some $ s\in G $ then $ g = h $. It follows from \cite{Riben} that every ordered monoid is cancellative.

\begin{lemma}\label{l2}
	Let $  R[G] $ be the monoid ring of a positive monoid $ G $ over a ring $ R $. If $ e = e_{0}g_{0} + e_{1}g_{1}+ \cdots + e_{m}g_{m} \in R[G] $ is an arbitrary idempotent satisfying $ g_{i } < g_{j} $ if $ i < j $, then $ g_{0} = \mu $, where $ \mu $ is the identity element of $ G $.
\end{lemma}

\begin{proof}
	It is known that $ e^{2} = \sum_{i}\sum_{j}e_{i}e_{j}g_{i}g_{j} $. So we have, $ e_{0}g_{0} = \sum_{i}\sum_{j} e_{i}e_{j}g_{i}g_{j} $. We show that, $ g_{0} = g^{2}_{0} $. Consider the following three cases:
	\begin{enumerate}
		\item If $ i = j = 0 $, then $ g_{0} = g^{2}_{0} $. 
		\item If $ i = 0 \neq j $, then $ i  < j $, and so $ g_{i} < g_{j} $. Thus, $ g^{2}_{0} = g_{i}g_{i} < g_{i}g_{j} = g_{0} $. On the other hand, $ g_{0} = g_{0}\mu < g_{0}g_{0} = g^{2}_{0} $. Therefore, $ g_{0} = g^{2}_{0} $.
		\item If $ i, j \neq 0 $, then $ 0 < i,j $. Hence, $ g_{0} < g_{i} $ and $ g_{0} < g_{j} $. Thus, $ g^{2}_{0} < g_{i}g_{j} = g_{0} $. On the other hand, $ g_{0} = g_{0}\mu < g_{0}g_{0} = g^{2}_{0} $. Therefore, $ g_{0} = g^{2}_{0} $.
	\end{enumerate}
Hence, $ g_{0}  = \mu $, since $ G $ is cancellative. Therefore, $ e_{0}\mu = e^{2}_{0}\mu $ and so $ e_{0} = e^{2}_{0} $.
\end{proof}

\begin{lemma}\label{l3}
	Let $  G $ be a positive monoid with identity element $ \mu $. Suppose $ e = e_{0}g_{0} + e_{1}g_{1}+ \cdots + e_{m}g_{m} $ is an idempotent element of $ R[G] $ satisfying $ g_{i } < g_{j} $ if $ i < j $. Then $ e_{i} \in Re_{0}R $ for all $  i \geq 0 $. 
\end{lemma}
\begin{proof}
	We prove the claim by induction on $ i $. For $ i = 0 $, we have $ e_{0} = e_{0}e_{0}$, by Lemma \ref{l2}. Hence, $ e_{0} \in Re_{0}R $. We assume the claim holds for $ 0\leq i < m-1 $. From Lemma \ref{l1} we have, $ e_{m}g_{m} = \sum_{i\leq m}\sum_{j\leq m}e_{i}e_{j}g_{i}g_{j}  $. Using the induction hypothesis, if at least one of $ i $ and $ j $ is less than $ m $, then $ e_{i}e_{j}\in Re_{0}R $. Moreover, since $ e_{m}e_{0}=e_{m}e_{0}e_{0}$ and $ e_{0}e_{m}=e_{0}e_{0}e_{m} $, we have $ e_{0}e_{m}, e_{m}e_{0}\in Re_{0}R $. On the other hand, the case $ i = j = m $ cannot occur because otherwise, $ g^{2}_{m} = g_{i}g_{j} = g_{m} $ and so $ g_{m} = \mu $, using cancellation property of $ G $.  Hence, $ e_{m}\in Re_{0}R $, and the result follows.
\end{proof}

\begin{theorem}\label{t1}
Let $ R[G] $ be the monoid ring of a positive monoid $ G $ over a ring $ R $. Then $ R $ is  right $\mathfrak{cP}$-Baer  if and only if $ R[G] $ is  right $\mathfrak{cP}$-Baer.
\end{theorem}
\begin{proof}
$ (\Rightarrow) $ Let $ e = e_{0}g_{0} +e_{1}g_{1}+\cdots+e_{m}g_{m}\in R[G] $ be an idempotent satisfying $ g_{0 } < g_{1}<\cdots < g_{m} $. It follows from Lemma \ref{l2}, $ e_{0} $ is an idempotent in $ R $. So there exists $ c = c^{2}\in R $ such that $ r_{R}(e_{0}R) = cR $. We shall prove $ r_{R[G]}(eR[G]) = cR[G] $. Let $ \alpha = a_{0}h_{0}+a_{1}h_{1}+\cdots+a_{n}h_{n} $ be an arbitrary element in $ r_{R[G]}(eR[G]) $ satisfying $ h_{0 } < h_{1}<\cdots < h_{n} $. Thus, $ er\alpha = 0 $ for all $ r\in R $. We proceed by induction on $ j $ to show that $ a_{j}\in r_{R}(e_{0}R) $ for each $ 0\leq j \leq n $. We have, $ 0 = er\alpha = e_{0}ra_{0}g_{0}h_{0} + \cdots + e_{m}ra_{n}g_{m}h_{n} $. Since $ g_{0}h_{0} $ is the smallest element in the $ g_{i}h_{j} $'s, its coefficient equals zero, that is,
$ e_{0}ra_{0} = 0 $. Thus, $ a_{0}\in r_{R}(e_{0}R) $ and so the result is true for $ j = 0 $. Now assume the result is true for all $ 0\leq j < n $. Hence, $ a_{j}\in r_{R}(e_{0}R) $ for any $ 0\leq j < n $. Thus, by Lemma \ref{l3}, $ a_{j}\in r_{R}(e_{i}R) $ for each $ 0\leq j < n $ and $ 0\leq i \leq m $. Now we have,
\begin{equation*}
	e_{0}ra_{n }+(\sum_{0\neq i\leq n}\sum_{j\leq n}e_{i}ra_{j}) = 0.
\end{equation*}
 In the above summation, the case $ j= n $ cannot occur because otherwise, $ g_{i}h_{n} = h_{n} $ and so $ g_{i} = \mu = g_{0} $, using cancellation property of $ G $, a contradiction. Now by the induction hypothesis, we have $ \sum_{0\neq i\leq n}\sum_{j< n}e_{i}ra_{j} = 0  $. Thus, $ e_{0}ra_{n } = 0 $, and so $ a_{n}\in r_{R}(e_{0}R) $. Then $ a_{j}\in cR $ and so $ ca_{j} = a_{j} $ for all $ 0\leq j\leq n $. Hence, $ c\alpha = c(\sum^{n}_{j=0}a_{j}h_{j}) = \sum^{n}_{j=0}ca_{j}h_{j} = \alpha $. Thus, $ \alpha\in cR[G] $ and so $ r_{R[G]}(eR[G]) \subseteq cR[G] $. Next, we claim that $ cR[G] \subseteq r_{R[G]}(eR[G]) $. It is sufficient to show $ c\in r_{R[G]}(eR[G]) $. Using Lemma \ref{l3}, we obtain $ c\in r_{R[G]}(eR[G]) $. Therefore, $R[G] $ is a right $\mathfrak{cP}$-Baer ring.\\  
 $ (\Leftarrow) $ Let $ R[G] $ is a right $\mathfrak{cP}$-Baer ring and $ c = c^{2}\in R $. Thus, $ c = c^{2} \in R[G] $, and so $ r_{R[G]}(cR[G]) = eR[G] $ where $ e = e_{0}g_{0} +e_{1}g_{1}+\cdots+e_{m}g_{m}\in R[G, <] $ is an idempotent $ g_{0 } < g_{1}<\cdots < g_{m} $. By Lemma \ref{l2}, it is enough to show that $ r_{R}(cR) = e_{0}R $. Let $ a\in r_{R}(cR)\subseteq r_{R[G]}(cR[G]) $. Thus, $ a\mu \in eR[G] $ and so $ ea\mu = a\mu $. Hence, $ e_{0}a = a $. Therefore, $ a\in e_{0}R $ and so $ r_{R}(cR)\subseteq e_{0}R $. Conversely, we have $ r_{R[G]}(cR[G]) = eR[G] $ so $ e\in r_{R[G]}(cR[G]) $. Thus, $ cRe = 0 $ and so $cRe_{0} = 0 $. Therefore, $ e_{0}R\subseteq r_{R}(cR) $, and that $ R $ is a right $\mathfrak{cP}$-Baer ring.
\end{proof}

It is well-known that torsion-free nilpotent groups and free groups are ordered groups (see \cite[Lemma 13.1.6 and Corollary 13.2.8]{pass}). Hence, we have the following. 

\begin{corollary}\label{c1}
Let $ G $ be a positive submonoid of a free or torsion-free nilpotent group. Then the monoid ring $ R[G] $ is a right $\mathfrak{cP}$-Baer ring if and only if $ R $ is a right $\mathfrak{cP}$-Baer ring.
\end{corollary}
 
 \begin{corollary}\label{c2}
 Let $ G $ be a submonoid of $ (N \cup {0})^{n} $ $ (n \geq 2) $, endowed with the order $ \leq $ induced by the product order, or lexicographic order or reverse lexicographic order. Then the monoid ring $ R[G] $ is a right $\mathfrak{cP}$-Baer ring if and only if $ R $ is a right $\mathfrak{cP}$-Baer ring.
 \end{corollary}

\begin{corollary}\label{c3}
	Let $ R $ be an arbitrary ring, $ G $ be either $ \mathbb{Q}^{+} = \left\lbrace a \in \mathbb{Q}\ |\ a \geq 0\right\rbrace  $ or $  \mathbb{R}^{+} = \left\lbrace a \in \mathbb{R}\ |\ a \geq 0\right\rbrace $ where $ \leq $ is the usual order. Then the monoid ring $ R[G] $ is a right $\mathfrak{cP}$-Baer ring if and only if $ R $ is a right $\mathfrak{cP}$-Baer ring.
\end{corollary}

\subsection{$\mathfrak{cP}$-Baer skew Laurent polynomial  rings }\label{s22}

Recall that, a ring $ R $ is semicommutative if $ ab = 0 $ implies $ aRb = 0 $ for $ a, b \in R $ (this ring is also called a zero insertion(ZI) ring in \cite{Habeb}, \cite {Kim} and \cite{WANG}). Semicommutative rings have also been called $IFP$ rings in the literature.  The ring $R$ is semicommutative if and only if any right (left) annihilator over $R$ is an ideal of R by \cite[Lemma 1]{HUH} or \cite[Lemma 1.2] {SHIN}. Every commutative ring is semicommutative.\par
According to Cohn \cite{cohn}, a ring $R$ is called reversible if $ab = 0$ implies $ba = 0,$ for $a, b \in R$. Before Cohn’s work, reversible rings were studied under the name completely reflexive by Mason in \cite{Mason}. Reduced rings (i.e., rings with no non-zero nilpotent elements) are reversible, and reversible rings are also semicommutative.\par
If $\alpha $ is a ring monomorphism, then the set $\{x^j \}_{j\geq0}$ is easily seen to be a left Ore subset of $R[x; \alpha]$. The skew Laurent polynomial ring  $R[x, x^{-1}; \alpha]$  is  the left quotient ring of $R[x; \alpha]$ with respect to $\{x^j\}_{j\geq0}$.  Elements of $R[x,x^{-1}; \alpha]$  are finite sums of elements of the form $x ^{-j} rx^i$, where $r \in R$ and $ i, j$ are nonnegative integers. Multiplication is subject to $xr = \alpha(r)x$ and $rx^{-1} = x^{-1}\alpha(r)$ for all $r \in R$.
\begin{theorem}\label{28}
	Let R be a semicommutative ring and $\alpha$ a compatible epimorphism of R. Then $ R $ is right $\mathfrak{cP}$-Baer  if and only if $ R[x,x^{-1};\alpha] $ is a right $\mathfrak{cP}$-Baer ring.
\end{theorem}
\begin{proof}
	By \cite[Theorem $ 3.5 $]{hamidi} for each semicommutative and $ \alpha $-compatible ring $ R $ we have, $ I(R) = I(R[x,x^{-1};\alpha])$. Using this and a similar argument as employed in Theorem \ref{12}, the result follows.
\end{proof}

\begin{corollary}
For a semicommutative ring $ R $, $ R $ is right $\mathfrak{cP}$-Baer  if and only if $ R[\mathbb{Z}] $ is a right $\mathfrak{cP}$-Baer ring.
\end{corollary}
\begin{proof}
	Theorem \ref{28} yields the result as $ R[x,x^{-1}]\cong R[\mathbb{Z}] $.
\end{proof}
  
Now we consider Jordan's construction of the ring $ A(R,\alpha) $  (See \cite{Jordan}, for more details) when $ \alpha $ is a monomorphism of $R$. Let $ A(R,\alpha) $ or $ A $ be the subset $ \lbrace x^{-i}ax^{i}\ |\ a\in R,\ i\geq 0\rbrace $ of the skew Laurent polynomial ring $ R[x,x^{-1};\alpha] $. For each $ j\geq0 $, we have $ x^{-i}ax^{i} = x^{-(i+j)}\alpha^{j}(a)x^{(i+j)} $. It then follows that the set of all such elements forms a subring of $ R[x,x^{-1};\alpha] $ with 
\begin{equation*}
	x^{-i}ax^{i} + x^{-j}bx^{j} = x^{-(i+j)}(\alpha^{j}(a) + \alpha^{i}(b)) x^{(i+j)},
\end{equation*}
and
\begin{equation*}
	(x^{-i}ax^{i})(x^{-j}bx^{j}) = x^{-(i+j)}(\alpha^{j}(a) \alpha^{i}(b)) x^{(i+j)},
\end{equation*}
for $ a,b\in R $ and $ i,j\geq0 $. Note that $ \alpha $ is actually an automorphism of $ A $. We have $ R[x,x^{-1};\alpha]\cong A[x,x^{-1};\alpha] $, by way of an isomorphism which maps $ x^{-i}ax^{j} $ to $ \alpha^{-i}(a)x^{j-i} $. Moreover, from \cite[Lemma 2.1]{hamidi}, $ I(A) = \lbrace x^{-i}ex^{i}\ |\ e\in I(R),\ i\geq 0\rbrace $.

\begin{lemma}\label{29}
	Let $ \alpha $ be a monomorphism of a ring $ R $ such that $ \alpha(e) = e $ for any left semicentral idempotent $ e\in R $. Then $ e\in S_{l}(R) $ if and only if $ x^{-i}ex^{i}\in S_{l}(A) $ for all $ i\geq 0 $.
\end{lemma}
\begin{proof}
	Let $ e\in S_{l}(R) $ and $ x^{-j}bx^{j}\in A $ be arbitrary. Thus, we have
	\begin{align*}
		(x^{-i}ex^{i})(x^{-j}bx^{j})(x^{-i}ex^{i}) &= (x^{-i}ex^{i})(x^{-(j+i)}\alpha^{i}(b)ex^{(j+i)})\\
		&= x^{-(i+j+i)}e\alpha^{i}(\alpha^{i}(b))ex^{(i+j+i)}\\
		&= x^{-(i+j+i)}\alpha^{i}(\alpha^{i}(b))ex^{(i+j+i)}\\
		&= x^{-(i+j+i)}\alpha^{i}(\alpha^{i}(b))\alpha^{i}(e)x^{(i+j+i)}\\
		&= x^{-(i+j)}\alpha^{i}(b)ex^{(i+j)}\\
		&= (x^{-j}bx^{j})(x^{-i}ex^{i}).
	\end{align*}
	Conversely, let  $ x^{-i}ex^{i}\in S_{l}(A) $. Thus, $ A(x^{-i}ex^{i}) = (x^{-i}ex^{i})A(x^{-i}ex^{i}) $ and $ r(x^{-i}ex^{i}) $\\$= (x^{-i}ex^{i})r(x^{-i}ex^{i}) $ for $ r\in R $. Then, $ x^{-i}\alpha^{i}(r)ex^{i} = x^{-(i+i)}\alpha^{i}(e)\alpha^{i}(\alpha^{i}(r))\alpha^{i}(e)x^{(i+i)} $. Hence, $ x^{-i}\alpha^{i}(r)ex^{i} = x^{-i}e\alpha^{i}(r)ex^{i} $ and so $ \alpha^{i}(r)e = e\alpha^{i}(r)e $. Since $ \alpha $ is a monomorphism, $ re = ere $ and we are done.
\end{proof}

\begin{theorem}\label{30}
	Let $ R $ be a ring with a compatible endomorphism $ \alpha $. Then R is a right $\mathfrak{cP}$-Baer ring if and only if $ A = A(R,\alpha) $ is a right $\mathfrak{cP}$-Baer ring.
\end{theorem}
\begin{proof}
	$ (\Rightarrow) $ Let R be right $\mathfrak{cP}$-Baer ring and $ x^{-i}ex^{i}\in A$ be an idempotent. Thus, $ e\in I(R) $ and, by definition, there exists $ c = c^{2}\in R $ such that $ r_{R}(eR) = cR $. As $ r_{R}(eR) $ is an ideal of $ R $, $ c\in S_{l}(R)  $ and so $ x^{-i}cx^{i}\in S_{l}(A) $, by Lemma \ref{29}. We show that $ r_{A}((x^{-i}ex^{i})A) = (x^{-i}cx^{i})A $. Let $ x^{j}bx^{j}\in r_{A}((x^{-i}ex^{i})A) $ and $ x^{-k}ax^{k}\in A $ be arbitrary. Therefore, $ (x^{-i}ex^{i})(x^{-k}ax^{k})(x^{j}bx^{j}) = 0 $ and so $ x^{-(i+k+j)}e\alpha^{i}(\alpha^{j}(a)\alpha^{k}(b))x^{(i+k+j)} = 0 $. Hence, $ e\alpha^{i}(\alpha^{j}(a)\alpha^{k}(b)) = 0 $, and so $ e\alpha^{j}(a)\alpha^{k}(b) = 0 $. As $ R $ is $ \alpha $-compatible, $ e\alpha^{j}(a\alpha^{k}(b)) = 0 $. Thus, $ ea\alpha^{k}(b) = 0 $, and so $ \alpha^{k}(b)\in r_{R}(eR) $. Hence, $  \alpha^{k}(b)\in cR $ and so $ c\alpha^{k}(b) = \alpha^{k}(b) $. Thus, $ c\alpha^{i}(b) = \alpha^{i}(b) $, since $ k\geq 0 $ is an arbitrary integer. Then we have,
	\begin{align*}
		(x^{-i}cx^{i})(x^{-j}bx^{j}) &= x^{-(i+j)}(\alpha^{j}(c)\alpha^{i}(b))x^{(i+j)})\\
		&= x^{-(i+j)}c\alpha^{i}(b)x^{(i+j)}\\
		&= x^{-(j+i)}\alpha^{i}(b)x^{j+i} = x^{-j}bx^{j}.
	\end{align*}
	Therefore, $ x^{-j}bx^{j}\in (x^{-i}cx^{i})A $ and so $r_{A}((x^{-i}ex^{i})A)\subseteq (x^{-i}cx^{i})A $. For the converse, it suffices to prove that $ x^{-i}cx^{i}\in r_{A}((x^{-i}ex^{i})A) $. Let $ x^{-j}bx^{j} $ is an arbitrary element of $ A $. Then $ (x^{-i}ex^{i})(x^{-j}bx^{j})(x^{-i}cx^{i}) = x^{-(i+j+i)}(e\alpha^{i}(\alpha^{i}(b))c) x^{(i+j+i)}$. Thus, $ x^{-i}cx^{i}\in r_{A}((x^{-i}ex^{i})A) $, since $ e\alpha^{i}(\alpha^{i}(b))c = 0 $. Therefore, $ (x^{-i}cx^{i})A\subseteq r_{A}((x^{-i}ex^{i})A) $ and so $ A $ is a right $\mathfrak{cP}$-Baer ring.\\
	$ (\Leftarrow) $ Let $ A $ is right $\mathfrak{cP}$-Baer and $ e = e^{2}\in R $. Thus, $ r_{A}(eA) = x^{-i}cx^{i}A$ for some idempotent element $ x^{-i}cx^{i}\in A $. We claim that $ r_{R}(eR) = cR $. As $ x^{-i}cx^{i}\in r_{A}(eA)$, $ erx^{-i}cx^{i} = 0 $ for each $ r\in R $. Thus, $ e\alpha^{i}(r)c = 0 $ and so $ erc = 0 $, by compatibility condition of $ \alpha $. Therefore, $ c\in r_{R}(eR) $ and so $ cR\subseteq r_{R}(eR) $. Conversely, assume that $ a\in r_{R}(eR) $ and $ x^{-j}rx^{j}\in A $. Thus, $ ex^{-j}rx^{j}a = x^{-j}er\alpha^{j}(a)x^{j} $. Since $ R $ is $ \alpha $-compatible and $ era = 0 $, $ er\alpha^{j}(a) = 0 $. Hence, $ a\in r_{A}(eA) $ and so $ a\in x^{-i}cx^{i}A $. Thus, $ a = ca $. Therefore, we are done with $ r_{R}(eR) \subseteq cR $.
\end{proof}

\begin{theorem}\label{31}
	Let R be a semicommutative ring and $ \alpha $ a compatible endomorphism of $ R $. Then $ R $ is a right $\mathfrak{cP}$-Baer ring if and only if $ R[x,x^{-1};\alpha] $ is a right $\mathfrak{cP}$-Baer ring.
\end{theorem}
\begin{proof}
	By using \cite[Lemma $2.3$]{hamidi}, $ A = A(R,\alpha) $ is a semicommutative ring. Since $ \alpha $ is an automorphism of $ A $ and $ R[x,x^{-1};\alpha]\cong A[x,x^{-1};\alpha] $, so the result follows by Theorems \ref{28} and \ref{30}.
\end{proof}
Using \cite [Corollary 2.5]{Alhevaz1}, as we have mentioned in Theorem \ref{321}, we can deduce that when $\alpha$ and $\sigma$ are rigid endomorphisms of a ring $R$ such that $\alpha\sigma = \sigma\alpha$, then the triangular matrix ring $T(R, n, \sigma)$ is a reversible and hence semicommutative $ \alpha$-compatible $\mathfrak{cP}$-Baer ring. We also note that any domain with an injective endomorphism $\alpha$ is an $\alpha$-rigid ring.\\
Now we provide various examples of $(\alpha,\delta)$-compatible  $\mathfrak{cP}$-Baer rings:\\
Let $R$ be a ring and $\sigma$ denotes an endomorphism of $R$ with $\sigma(1) = 1$. Let
\begin{equation*}
	T_{n}(R,\sigma) = \lbrace (a_{ij})_{n\times n}\ |\  a_{ij}\in R,\ a_{ij} = 0\ if\ i>j\rbrace .
\end{equation*}
For $ (a_{ij}), (b_{ij})\in T_{n}(R,\sigma) $ define,
\begin{equation*}
	(a_{ij}) + (b_{ij}) = (a_{ij}+b_{ij})\quad and\quad (a_{ij}) \ast (b_{ij}) = (c_{ij}),
\end{equation*}
where $ c_{ij} = 0 $ for $ i>j $ and $ c_{ij} = \sum^{j}_{k=i}a_{ik}\sigma^{k-i}(b_{kj}) $ for $ i\leq j $.  It can be easily checked that $ T_{n}(R, \sigma) $ is a ring, called the skew triangular matrix ring over R. If $ \sigma $ is the identity, $ T_{n}(R,\sigma) $ is the triangular matrix ring  $ T_{n}(R) $.\\
The subring of the skew triangular matrices  with constant main diagonal is denoted by $S(R,n,\sigma)$; and the subring of the skew triangular matrices with constant diagonals is denoted by $T(R,n,\sigma)$. We can denote $A=(a_{ij})\in T(R,n,\sigma)$ by $(a_{11},\ldots,a_{1n})$. Then $T(R,n,\sigma)$ is a ring with addition point-wise and multiplication given by:
\begin{equation*}
	(a_{0},\ldots, a_{n-1})(b_{0},\ldots, b_{n-1})= (a_{0}b_{0}, a_{0}*b_{1}+a_{1}*b_{0},\ldots, a_{0}*b_{n-1}+\cdots +a_{n-1}*b_{0}),
\end{equation*}
 with $a_{i}*b_{j} = a_{i}\sigma^{i}(b_{j})$, for each $i$ and $j$. Also, there exists a ring isomorphism $$  \varphi: R[x;\sigma]/ \langle x^{n}\rangle \longrightarrow T(R,n,\sigma) $$ such that $ \varphi(a_{0}+a_{1}x+\cdots+a_{n-1}x^{n-1}+\langle x^{n}\rangle) = (a_{0},\ldots, a_{n-1}) $. Therefore, 
 \begin{equation*}
 	T(R,n,\sigma)\cong R[x;\sigma]/\langle x^{n}\rangle ,
 \end{equation*}
 where $\langle x^{n}\rangle$ is the ideal generated by $x^{n}$ in $R[x;\sigma]$.\\
We consider the following two subrings of $S(R,n,\sigma)$, as follows:  (see\cite{Zhou1})
\begin{equation*}
	A(R,n,\sigma) = \sum_{j=1}^{\lfloor \frac{n}{2}\rfloor}\sum_{i=1}^{n-j+1}{a_{j}E_{i,i+j-1}} + \sum_{j=\lfloor \frac{n}{2}\rfloor+1}^{n}\sum_{i=1}^{n-j+1}{a_{i,i+j-1}E_{i,i+j-1}};
\end{equation*}
\begin{equation*}
	B(R,n,\sigma)=\{A+rE_{1k}\,|\,A\in A(R,n,\sigma)\,\, \mathrm{and}\,\,r\in R\}\qquad where,\ n=2k\geq 4.
\end{equation*}

Let $\alpha$ and $\sigma$ be endomorphisms of a ring $R$ and $\delta$ is an $\alpha$-derivation, with $\alpha\sigma=\sigma\alpha$ and $\delta\sigma=\sigma\delta$. The endomorphism $\alpha$ of $R$ is extended to the endomorphism $\bar{\alpha} : T_n(R,\sigma) \rightarrow T_n(R,\sigma)$ defined by $\bar{\alpha}((a_{ij})) = (\alpha(a_{ij}))$ and the $\alpha$-derivation $\delta$ of $R$ is also extended to $\bar{\delta}: T_n(R,\sigma)\rightarrow T_n(R,\sigma)$ defined by $\bar{\delta}((a_{ij})) = (\delta(a_{ij}))$.

\begin{proposition}\label{32}
	Let $R$ be a domain  with a  monomorphism $ \sigma $. Then  the skew matrix rings $T(R,n,\sigma)$, $A(R,n,\sigma)$ and $ B(R,n,\sigma) $ are right $\mathfrak{cP}$-Baer rings.
\end{proposition}
\begin{proof}
	By \cite[Theorem $ 2.8 $]{Alhevaz1},  for the $\sigma$-rigid ring $R$, the matrix rings  $A(R,n,\sigma)$, $B(R,n,\sigma)$
	and $T(R,n,\sigma)$ are Armendariz rings. Using \cite[Lemma $ 7 $]{Kim}, Armendariz rings are abelian, so the proof follows.
\end{proof}

\begin{proposition}\label{321}
 Let $\sigma$ and $\alpha$ be rigid endomorphisms and $\delta$ an $\alpha$-derivation of a ring $R$ such that $\alpha\sigma=\sigma\alpha$ and $\delta\sigma = \sigma\delta$. Then $A(R,n,\sigma)$, $B(R,n,\sigma)$ and $T(R,n,\sigma)$ are $(\bar{\alpha},\bar{\delta})$-compatible right $\mathfrak{cP}$-Baer rings.
\end{proposition}

\begin{proof}
It follows from Proposition \ref{32} and \cite[Theorem $ 2.3 $]{Alhevaz1}.
\end{proof}

From \cite{b1}, we see there are close relationships between quasi-Baer rings and $FI$-extending rings. Also, in \cite{birkenkp}, the authors found a relationship between p.q.-Baer rings and $pFI$-extending rings. Recall that $R$ is $FI$-extending ($pFI$-extending) if every ideal (principal ideal) is essential in a direct summand of $R_R$. Naturally, one may ask: Is there a corresponding extending condition that is in close relationship to the right $\mathfrak{cP}$-Baer condition? Since $\mathfrak{cP}$-Baer rings are defined in terms of principal ideals generated by idempotents, our extending condition should be similarly defined. By \cite[Definition $ 1.9 $]{Heider},  a ring $R$ is called right $I$-extending if every principal ideal of $R$ that is generated by an idempotent is essential in a direct summand of $ R_{R} $.

\begin{definition}\cite[Definition $ 1.9 $]{Heider}\label{14}
	A ring $R$ is said to be right (left) $ I $-$\mathit{extending}$ if for each $ e = e^{2}\in R $ there exists idempotent $ c\in R $ such that $ ReR_{R}\leq^{ess}cR_{R} $ $ (respectively,\ _{R}ReR\leq^{ess}\ _{R}Rc) $. Also, $R$ is called $ I $-extending if it is both left and right $ I $-extending.
\end{definition}
In \cite{Heider}, the authors give basic results for $\mathfrak{cP}$-Baer rings, $I$-extending rings, and conditions which ensure that the $\mathfrak{cP}$-Baer and $I$-extending properties are equivalent.
By \cite[Theorem $ 1.11 $]{Heider}, when  $R$ is a semiprime ring, then $R$ is right $\mathfrak{cP}$-Baer if and only if $R$ is right $ I $-extending. Also, by \cite [Corollary $ 1.12 $]{Heider}, when $R$ is a semiprime ring, then $R$ is right $ I $-extending if and only if $R$ is left $ I $-extending.\\
Using the above notes, we have the following.
\begin{lemma}\label{ll1}
Let $ R $ be a semiprime ring with a compatible endomorphism $ \alpha $. Assume $ f(x) = \sum_{i = 0}^{m}a_{i}x^{i} $ and $ g(x) = \sum_{j = 0}^{n}b_{j}x^{j} $ are arbitrary elements of $ R[x;\alpha] $. Then $ f(x)R[x;\alpha]g(x) = 0 $ if and only if $ a_{i}Rb_{j} = 0 $ for each $ i, j $. 
\end{lemma}
\begin{proof}
$ (\Rightarrow) $ If $ f(x)R[x;\alpha]g(x) = 0 $, then $ \sum_{i = 0}^{m}\sum_{j = 0}^{n}a_{i}\alpha^{i}(rb_{j})x^{i+j} = 0 $ for each $ r \in R $. Hence, $ a_{m}\alpha^{m}(rb_{n}) = 0 $, and so $ a_{m}Rb_{n} = 0 $, as $ R $ is $ \alpha$-compatible. Now, replace $ r $ with $ rb_{n-1}sa_{m}t $, where $ r, s, t $ are elements of $ R $. Thus we have $ \sum_{i = 0}^{m}\sum_{j = 0}^{n-1}a_{i}\alpha^{i}(rb_{n-1}sa_{m}tb_{j})x^{i+j} = 0 $. Therefore, $ a_{m}\alpha^{m}(rb_{n-1}sa_{m}tb_{n-1}) = 0 $ and hence,  $ a_{m}rb_{n-1}sa_{m}tb_{n-1} = 0 $. Thus, $ (Ra_{m}Rb_{n-1})^{2} = 0 $. As $ R $ is semiprime, $ a_{m}Rb_{n-1} = 0 $. By a similar argument, we obtain $ a_{m}Rb_{j} = 0$, for $ 0\leq j\leq m $. From $ \alpha $-compatibility of $ R $, we get $ \sum_{i = 0}^{m-1}\sum_{j = 0}^{n}a_{i}\alpha^{i}(rb_{j})x^{i+j} = 0 $. Using induction on $ m + n $, we get $ a_{i}Rb_{j} = 0 $ for each $ i, j $. \\
$ (\Leftarrow) $ It follows by a similar argument.
\end{proof}
\begin{proposition}\label{17}
	Let $ R $ be a  ring and $ \alpha $ be a compatible endomorphism of $R$. The following statements are equivalent:
	\begin{itemize}
	\item[(i)] $R$ is semiprime.
	\item[(ii)] $ R[x;\alpha] $ is semiprime.
	\item[(iii)] $ R[x,x^{-1};\alpha] $ is semiprime.
	\end{itemize}
\end{proposition}

\begin{proof}
(i) $ \Rightarrow $ (ii) Let $ f(x) = \sum_{i= 0}^{m}a_{i}x^{i} $ be a non-zero element of $ R[x;\alpha] $ such that $ f(x)R[x;\alpha]f(x) = 0 $. Using Lemma \ref{ll1}, $ a_{i}Ra_{i} = 0 $, for any $ i\geq 0 $. Since $ R $ is semiprime, $ a_{i} = 0 $, for each $ i\geq 0 $. Therefore, $ f(x) = 0 $ and so $ R[x;\alpha] $ is semiprime.\\
(ii) $ \Rightarrow $ (iii) Assume that $ R[x,x^{-1};\alpha] $ is not  semiprime. Hence, $ \mathcal{P}(R[x,x^{-1};\alpha]) \neq 0 $. Let  $ I $ be a non-zero nilpotent ideal of $ R[x,x^{-1};\alpha] $. Hence, $ I\cap R[x;\alpha] $ is a non-zero nilpotent ideal of $ R[x;\alpha] $. Therefore, $ \mathcal{P}(R[x;\alpha]) \neq 0 $, a contradiction.\\
(iii) $ \Rightarrow $ (i) Suppose that $R$ is not semiprime. Thus, $ aRa = 0 $ for some non-zero element $ a\in R $. Using  $ \alpha $-compatibility of $R$, it implies that $ aR[x,x^{-1};\alpha]a = 0 $, which is a contradiction.
\end{proof}

\begin{corollary}\label{12o}
	Let R be an $ \alpha $-compatible semiprime ring with an endomorphism $ \alpha $. Then the following statements are equivalent: 
	\begin{enumerate}
		\item R is $ I $-extending;
		\item $ R[x;\alpha] $ is $ I $-extending;
		\item $ R[x,x^{-1};\alpha] $ is $ I $-extending.
			\end{enumerate}
\end{corollary}

\subsection{$\mathfrak{cP}$-Baer  skew inverse power series rings}\label{s3}

We denote by $R((x^{-1}; \alpha, \delta))$ the skew inverse  Laurent series ring over the coefficient ring $R$ formed by formal series $f(x)=\sum^{\infty}_{ i= m} a_ix^{-i}$ where $x$ is a variable, $m$ is an integer (maybe negative), and the coefficients $a_i$ of the series $f(x)$ are elements of the ring $R$. In the ring $R((x^{-1}; \alpha, \delta))$, addition is defined as usual and multiplication is defined with respect to the relations $xa =\alpha(a)x +\delta(a), x^{-1}a=\sum^{\infty}_{ i=1}  \alpha^{-1}(-\delta\alpha^{-1})^{i-1}(a)x^{-i}$,  for each $a \in R. $ Observe that the subset $R[[x^{-1}; \alpha, \delta]]$ of $R((x^{-1}; \alpha, \delta))$ consists of inverse power series of the form $f(x)=\sum^{\infty}_{ i=0} a_ix^{-i}$.\\
The interest in this notion derived initially from the fact that skew inverse Laurent series rings have wide applications.  The ring-theoretical properties of  skew  inverse Laurent series rings have been investigated by many authors (see \cite{Goodearl-pseudo}, \cite {Letzter}, \cite {Paykan.CMJ}, \cite{Psemiprime}, \cite{paykan1}, \cite{paykan}, \cite{tug}, \cite{tug-alpha}, \cite{tug2},  for instance). \\
By \cite[Examples 4.1. and 4.2]{Habibi}, there are examples which show that principally quasi Baerness of a ring $R$ and the skew inverse Laurent series ring $R((x^{-1}; \alpha, \delta))$ are not dependent, in general.

\begin{lemma}\label{19}
	Let R be an $(\alpha,\delta)$-compatible ring with an automorphism $ \alpha $ and an $ \alpha $-derivation $ \delta.$  If $ e(x) = \sum^{\infty}_{i = 0}e_{i}x^{-i} $ is an idempotent in $ R[[x^{-1}; \alpha, \delta]] $ then $ e_{i}\in Re_{0}R $ for any $ i\geq 0. $
\end{lemma}
\begin{proof}
	The proof is similar to Lemma \ref{2}.
\end{proof}

\begin{proposition}\label{20}
	Let R be an $(\alpha,\delta)$-compatible ring with an automorphism $ \alpha $ and an $ \alpha $-derivation $ \delta.$ If $ e(x)\in S_{l}(R[[x^{-1}; \alpha, \delta]]) $ where $ e(x) = \sum^{\infty}_{i = 0}e_{i}x^{-i}, $ then:
	\begin{enumerate}
		\item $ e_{0}\in S_{l}(R) $,
		\item $ e_{0}e_{i} = e_{i} $ for all $ i\geq 0 $,
		\item $ e_{i}e_{0} = 0 $ for all $ i\geq 1 $,
		\item $ e(x)R[[x^{-1}; \alpha, \delta]] = e_{0}R[[x^{-1}; \alpha, \delta]] $.
	\end{enumerate}
\end{proposition}
\begin{proof}
	The proof is similar to Proposition \ref{4}.
\end{proof}

\begin{lemma}\label{21}
	Let R be a right $\mathfrak{cP}$-Baer ring with an automorphism $ \alpha $ and an $ \alpha $-derivation $ \delta $ such that $ \alpha (c) = c $ for any left semicentral idempotent $ c\in R $. If $ eR\alpha^{k}(a) = 0 $ for some integer $ k $, where $ e = e^{2}\in R $ and $ a\in R $, then $ eR\alpha^{m_{1}}\delta^{n_{1}} \cdots \alpha^{m_{t}}\delta^{n_{t}}(a) = 0 $ for any integer $ m_{i}$ and any positive integer $ n_{j} $, $1\leq i, j\leq t$.
\end{lemma}
\begin{proof}
	We first show that $ eR\alpha^{m}(a) = 0 $ for any integer $ m $. Let $ eR\alpha^{k}(a) = 0 $ for some integer $ k $, where $ e = e^{2}\in R $ and $ a\in R $. Thus, by assumption there exists $ c = c^{2} \in R $ such that $ r_{R}(eR) = cR $. So $ c\in S_{l}(R) $ and $ \alpha^{k}(a)\in cR = \alpha^{k}(cR) $. Hence, $ a\in cR $ and so $ a = ca $. Therefore, $ eR\alpha^{m}(a) = eR\alpha^{m}(ca) = eRc\alpha^{m}(a) = 0 $ for any integer $ m $. We will show that $ eR\delta^{n}(a) = 0 $ for  any positive integer $ n $. By the preceding argument,  $ a = ca $. On the other hand, $ \delta(c) = c\delta(c) + \delta(c)c $. So $ \delta(c)c = c\delta(c)c + \delta(c)c $. As $ c $ is a left semicentral idempotent, $ c\delta(c)c = \delta(c)c $ and so $ \delta(c)c = 0 $. Thus, $ \delta(c) = c\delta(c) $. Therefore, $ eR\delta(a) = eR\delta(ca) = eR(c\delta(a) + c\delta(c)a) = 0 $. By induction, we can show that $ eR\delta^{n}(a) = 0 $ for any positive integer
	$ n $. As $ eR\alpha^{m}(a) = 0, $  $ eR\delta^{n}\alpha^{m}(a) = 0 $ for any positive integer $ n $. On the other hand, $ eR\delta^{n}(a) = 0 $ thus $ \delta^{n}(a)\in cR $ and so $ \delta^{n}(a) = c\delta^{n}(a) $. Hence, $ eR\alpha^{v}(\delta^{n}(a)) = eR\alpha^{v}(c\delta^{n}(a)) = eRc\alpha^{v}(\delta^{n}(a)) = 0 $ for any integer $ v $. By continuing this process, we have $ eR\alpha^{m_{1}}\delta^{n_{1}} \cdots \alpha^{m_{t}}\delta^{n_{t}}(a) = 0 $ for any integer $ m_{i}$ and any positive integer $ n_{j} $, $1\leq i, j\leq t$.
\end{proof}

We can prove the following using a similar argument as employed in  Proposition \ref{8}. 

\begin{lemma}\label{22}
	Let R be a right $\mathfrak{cP}$-Baer ring with an automorphism $ \alpha $ and an $ \alpha $-derivation $ \delta $ such that $ \alpha (c) = c $ for any left semicentral idempotent $ c\in R $. For any idempotent $ e(x) = \sum^{\infty}_{i=0}e_{i}x^{-i}\in  R[[x^{-1}; \alpha, \delta]]$ and any $ p(x) =\sum^{\infty}_{j=0}a_{j}x^{-j}\in  R[[x^{-1}; \alpha, \delta]] $, if $ e(x)R[[x^{-1}; \alpha, \delta]]p(x) = 0 $ then for any integer $ m_{l}$ and any positive integer $ n_{k} $, $1\leq l, k\leq t,  $ $e_{0}R\alpha^{m_{1}}\delta^{n_{1}} \cdots \alpha^{m_{t}}\delta^{n_{t}}(a_{j}) = 0 $ for each $ j\geq 0 $.
\end{lemma}

\begin{lemma}\label{23}
	Let R be a right $\mathfrak{cP}$-Baer ring with an automorphism $ \alpha $ and an $ \alpha $-derivation $ \delta $ such that $ \alpha (c) = c $ for any left (resp. right) semicentral idempotent $ c\in R $. If $ c $ is a left (resp. right) semicentral idempotent of R, then $ c $ is also a left (resp. right) semicentral idempotent in $ R[[x^{-1}; \alpha, \delta]] $.
\end{lemma}
\begin{proof}
	We will prove the case for left semicentral idempotent, and the other case is similar.
	As $ (c-1)Rc = 0, $ $ (c-1)R\alpha^{m_{1}}\delta^{n_{1}} \cdots \alpha^{m_{t}}\delta^{n_{t}}(c) = 0, $ for any integer $ m_{i}$ and any positive integer $ n_{j}, $ By Lemma \ref{21}, we have $ c\alpha^{-1}(\delta\alpha^{-1})^{d}(c) = 0 $ for each integer $ d>0 $. We will show that $ cx^{-n}c = x^{-n}c$ for any positive integer $ n $. If $ n= 1 $ then we have,
	\begin{equation*}
		(c-1)x^{-1}c = (c-1)(\alpha^{-1}(c)x^{-1} - \alpha^{-1}\delta\alpha^{-1}(c)x^{-2}+\cdots) 
		= 0.
	\end{equation*}
	If $ n=2, $ then
	\begin{align*}
		cx^{-2}c &=cx^{-1}x^{-1}c = cx^{-1}(cx^{-1}c)\\
		& = (cx^{-1}c)x^{-1}c = x^{-1}(cx^{-1}c)\\
		&= x^{-2}c.
	\end{align*}
	By continuing this way, we get $ cx^{-n}c = x^{-n}c$ for any positive integer $ n $. Now assume that $ p(x) =\sum^{\infty}_{j=0}a_{j}x^{-j}\in  R[[x^{-1}; \alpha, \delta]] $ is arbitrary. Hence,
	\begin{equation*}
		cp(x)c = \sum^{\infty}_{j=0}ca_{j}x^{-j}c = \sum^{\infty}_{j=0}ca_{j}cx^{-j}c = \sum^{\infty}_{j=0}a_{j}cx^{-j}c = p(x)c.
	\end{equation*}
	Therefore, $c$ is a left semicentral idempotent of $ R[[x^{-1}; \alpha, \delta]] $, and the result follows.
\end{proof}

\begin{theorem}\label{24}
	Let $R$ be an $ (\alpha,\delta) $-compatible ring  with an automorphism $ \alpha $ and an  $ \alpha $-derivation $ \delta $. Then R is right $\mathfrak{cP}$-Baer if and only if $ R[[x^{-1}; \alpha, \delta]] $ is right $\mathfrak{cP}$-Baer.
\end{theorem}
\begin{proof}
	Using Lemmas \ref{19}, \ref{20}, \ref{21}, \ref{22} and \ref{23}, the proof follows in a similar method as that employed in the proof of Theorem \ref{12}. 
\end{proof}

\begin{corollary}\label{24o}
	Let $R$ be an $ (\alpha,\delta) $-compatible abelian ring  with an automorphism $ \alpha $ and an  $ \alpha $-derivation $ \delta $. Then  $ R[[x^{-1}; \alpha, \delta]] $ is right $\mathfrak{cP}$-Baer.
\end{corollary}
Since each reduced ring is a $\delta$-compatible $\mathfrak{cP}$-Baer ring, we have

\begin{corollary}\label{24oo}
	Let $R$ be an $ \alpha$-compatible reduced ring  with an automorphism $ \alpha $ and an  $ \alpha $-derivation $ \delta $. Then  $ R[[x^{-1}; \alpha, \delta]] $ is right $\mathfrak{cP}$-Baer.
\end{corollary}
\begin{corollary}\label{24ooo}
	For any reduced ring  with a derivation $ \delta $,   $R[[x^{-1}; \delta]] $ is a $\mathfrak{cP}$-Baer ring.
\end{corollary}

\begin{proposition}\label{241}
	If $R$ is a $ \alpha $- compatible prime $ (resp., semiprime) $ ring, then  $R(({x^{ - 1}};\alpha ,\delta ))$ and $R[[{x^{ - 1}};\alpha ,\delta ]]$ are  prime $ (resp., semiprime) $ rings.
\end{proposition}
\begin{proof}
	Let $R$ be a prime ring and $f(x)=\sum^{\infty}_{ i= m} a_ix^{-i}$ and $g(x)=\sum^{\infty}_{ j=n} b_jx^{-j}$ be non-zero elements of $R(({x^{ - 1}};\alpha ,\delta ))$ such that $ a_{m}, b_{n}\neq 0 $. If $ f(x)R(({x^{ - 1}};\alpha ,\delta ))g(x) =0 $ then $ f(x)Rg(x) = 0 $. Thus, $ a_{m}x^{-m}Rb_{n} = 0 $ and so $ a_{m}Rb_{n} = 0 $ by compatibility condition, which contradicts our assumption, and we may assume that $ f(x)R(({x^{ - 1}};\alpha,\delta ))g(x) \neq 0 $ for all non-zero elements of $R(({x^{ - 1}};\alpha,\delta ))$. Therefore, $R(({x^{ - 1}};\alpha ,\delta ))$ is a prime ring. Similarly, we can show that $R(({x^{ - 1}};\alpha ,\delta ))$ is a semiprime ring.
\end{proof}	

\begin{corollary}\label{242}
	If $R$ is an $ \alpha $-compatible prime ring, then  $R(({x^{ - 1}};\alpha ,\delta ))$ and $R[[{x^{ - 1}};\alpha ,\delta ]]$ are  right $\mathfrak{cP}$-Baer rings.
\end{corollary}

\section*{Acknowledgement}
We would like to thank the referee for his/her carefully reading the paper and his/her helpful suggestions which improved its quality. 

\bibliographystyle{amsplain}

\end{document}